\documentclass[11pt]{article}

\usepackage[a4paper,margin=2.5cm]{geometry}

\usepackage[T1]{fontenc}
\usepackage[utf8]{inputenc}
\usepackage{lmodern}

\usepackage{amsmath,amssymb,amsthm}
\usepackage{graphicx}
\usepackage{booktabs}
\usepackage{multirow}
\usepackage{array}
\usepackage{enumitem}
\usepackage{float}
\usepackage{setspace}
\usepackage{caption}
\usepackage{subcaption}
\usepackage{hyperref}
\usepackage[nameinlink]{cleveref}
\usepackage{natbib}

\hypersetup{
    colorlinks=true,
    linkcolor=blue,
    citecolor=blue,
    urlcolor=blue
}

\title{\textbf{Question-Guided Learning for Mathematics Education: An Autonomous Online Inquiry-Based Approach}}

\author{
Borbála Neogrády-Kiss\thanks{Most of the research reported in this paper was conducted while the author was affiliated with Eötvös Loránd University.}\\
Budafok-Tétényi Baross Gábor Primary School\\
Budapest, Hungary
}

\date{}

\begin{document}

\maketitle

\begin{abstract}

This paper presents an autonomous online application designed to support inquiry-based learning in mathematics when in-person instruction is not possible. The application teaches the concept of boundedness of functions through adaptive questioning, requiring active problem solving without teacher intervention. It was tested with 92 Hungarian secondary school students preparing for the advanced-level final examination.

Assessment based on Bloom's taxonomy showed that most students achieved at least a medium level of understanding, with many reaching higher cognitive levels. The results suggest that such applications may support learning and complement traditional instruction.

Moreover, the study suggests that future AI-based educational systems may benefit from similar inquiry-oriented approaches, emphasizing active thinking and guided problem solving rather than passive answer generation.

\bigskip


\end{abstract}

\section{Introduction}

Supporting autonomous mathematics learning remains a significant challenge in both secondary and higher education. Although online learning environments provide flexibility, they often struggle to maintain the active engagement and independent thinking necessary for deeper mathematical understanding. Inquiry-based learning approaches may address this problem by encouraging students to actively explore questions, construct definitions, and solve problems rather than passively receive information. However, implementing such approaches in autonomous online environments remains difficult.

At the same time, recent developments in artificial intelligence have increased interest in AI-supported educational systems. While large language models can often generate convincing mathematical explanations, they may also generate mathematically incorrect or internally inconsistent reasoning. Moreover, current AI systems frequently emphasize answer generation rather than guiding students through the process of mathematical discovery and reflection. Consequently, an important question is whether educational technologies can be designed in ways that support active thinking through adaptive questioning and guided inquiry.

To explore this possibility, we developed an autonomous online learning application designed to teach the concept of boundedness of functions through adaptive questioning. The system guides students through progressively structured questions and adapts its responses depending on the types of errors students make, thereby encouraging active problem solving without teacher intervention. Rather than directly providing solutions, the application attempts to support conceptual understanding through guided inquiry and differentiated feedback.

The application was tested with 92 Hungarian secondary school students preparing for the advanced-level mathematics final examination. Students completed the learning material independently, after which their understanding was assessed using tasks aligned with the revised Bloom taxonomy. In addition to evaluating learning outcomes, we also examined the relationship between students' performance and the pace at which they progressed through the application.

To strengthen the evaluation of the application, learning outcomes were also compared with those of a matched control group completing the same assessment without prior instruction on boundedness.

The study addresses the following research questions:

\begin{enumerate}
    \item Does the developed application improve students' understanding of boundedness compared with students who do not receive instruction on the topic, and what levels of understanding can students achieve using the application?

    \item Is there a relationship between students' performance and the pace at which they progress through the application?
\end{enumerate}

\section{Theoretical Background}

Distance learning offers flexibility and may be particularly useful in situations where regular in-person participation is difficult, although reduced motivation and the lack of personal interaction may also create challenges \citep{Galusha1997}. At the same time, online learning environments can be effective: \citet{Soffer2018} found that students completing online courses achieved results comparable to, or in some cases better than, those participating in face-to-face instruction. Research has also suggested that teachers are generally able to adapt to online educational environments even without extensive prior ICT experience \citep{Spoel2020}.

Artificial intelligence has increasingly become part of educational practice \citep{Chen2020}. In this context, AI refers to computer systems capable of performing tasks associated with intelligent human behaviour \citep{RussellNorvig2010}. AI-supported educational systems may provide personalized feedback, identify students' learning difficulties, and support mathematical learning processes \citep{HwangTu2021}. Existing systems such as iTalk2Learn, ASSISTments, and AIDA illustrate how intelligent systems can support mathematics learning through adaptive feedback and guided assistance \citep{Grawemeyer2015,Heffernan2014,Chaudhry2021}. More broadly, the aims of AI in education include reducing teachers' workload, providing contextualized learning experiences, supporting students' metacognitive awareness, and creating intelligent learning environments capable of adaptive interaction \citep{Chaudhry2021,Ma2014}.

Contemporary AI models are increasingly capable of solving and explaining mathematical problems across different levels of complexity. The pace of development is illustrated by recent systems that have contributed to highly complex mathematical research, including work related to the disproof of the unit distance conjecture \citep{Alon2026,OpenAI2026}. At the same time, concerns have also emerged regarding the educational use of large language models. Recent research by \citet{Kosmyna2025} suggests that extensive reliance on AI-generated content may reduce cognitive engagement and weaken reflective learning processes. Although such findings should be interpreted cautiously, they nevertheless highlight the importance of designing AI-supported educational systems in ways that promote active thinking and meaningful student participation rather than passive dependence on automatically generated answers.

These concerns are closely connected to inquiry-based learning (IBL), which emphasizes active exploration, questioning, hypothesis testing, and student-centered knowledge construction \citep{Gillies2020}. Inquiry-based learning environments encourage students to investigate problems independently rather than passively receive information \citep{Artigue2013,Maass2013}. In particular, level-two inquiry, where students work with guided questions and structured support while actively constructing conclusions themselves, appears especially effective in mathematics education \citep{Bruder2013}. Research has shown that inquiry-oriented instruction may increase both engagement and mathematical performance compared to more traditional approaches \citep{Freeman2014,Kogan2014,Rasmussen2007}, and may also be more effective than other forms of inquiry instruction \citep{Jiang2015}.

However, successful inquiry-based learning requires carefully structured guidance and productive struggle \citep{Hiebert2007,Kapur2008,Russo2024}. Research involving pre-service teachers suggests that inquiry-oriented environments may initially create difficulties, although participants later tend to evaluate them positively and recognize their effectiveness \citep{Biber2024}. Teachers' own inquiry experiences are also important for effectively supporting students \citep{Makar2007}. These considerations are particularly relevant in computer-based learning environments, where guidance must be embedded into the structure of the application itself. In this regard, \citet{Hidayah2018} also emphasize the importance of manipulative-based teaching and guided question sequences.

The application developed in this study focuses on the concept of boundedness of functions. Previous research has shown that students often struggle to distinguish between maxima and upper bounds, minima and lower bounds, as well as between monotonicity and boundedness \citep{Bloch2003}. Difficulties also arise in understanding unboundedness and in constructing appropriate examples. The instructional design of the application was therefore partly based on these previously identified conceptual difficulties.

The research was also informed by the revised version of Bloom's taxonomy developed by \citet{Anderson2001}, which categorizes cognitive processes into six levels. This revised framework has previously been applied in educational technology research \citep{Crompton2019} and has been shown to be useful for evaluating levels of mathematical understanding \citep{Radmehr2019}. The present study therefore uses this framework to assess students' cognitive achievement and to compare learning outcomes between students using the application and a matched control group.

\begin{figure}[htbp]
    \centering
    \includegraphics[width=0.8\textwidth]{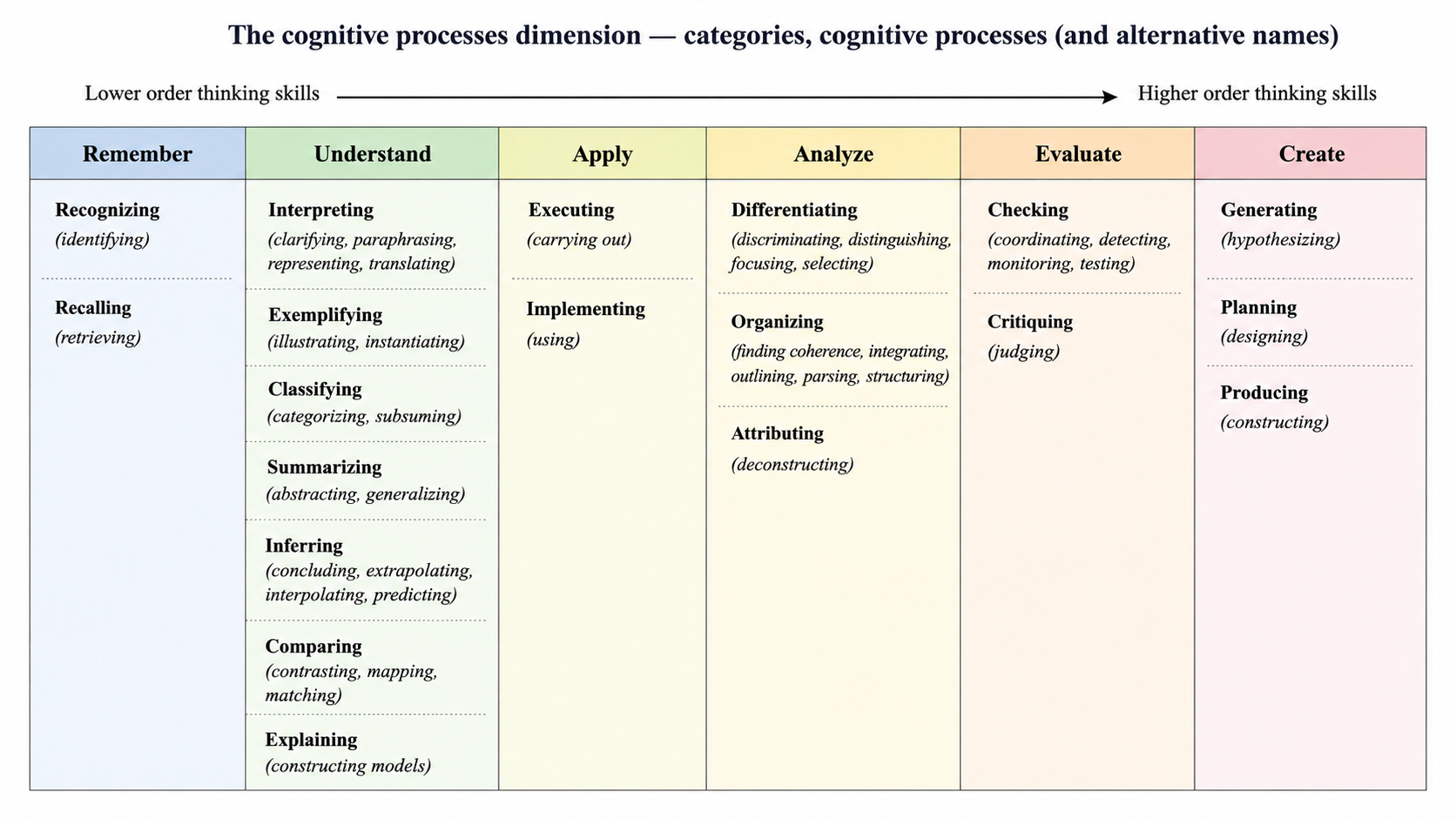}
    \caption{The cognitive process dimension of the revised Bloom taxonomy (adapted from Anderson et al., 2001).}
    \label{fig:bloom}
\end{figure}

The study also builds on our earlier pilot research involving pre-service mathematics teachers, where initial results suggested that the application may effectively support conceptual understanding of boundedness \citep{Neogrady2023}. While the pilot study provided encouraging preliminary evidence, it did not include a comparison group. The present study extends this earlier work by evaluating the application in a substantially larger sample and by comparing students' learning outcomes with those of a matched control group completing the same assessment without prior instruction.

\section{Structure of the Application}

The application was designed as an autonomous, inquiry-oriented learning environment for teaching the concept of boundedness of functions. Students progress through a sequence of adaptive questions without teacher intervention. Depending on their responses, the system either presents the next conceptual step or generates additional guiding questions intended to address the specific misconception reflected in the student's answer.

The revised version used in this study was developed for secondary school students aged 14--18. The application consists of 13 main tasks covering the central concepts related to boundedness, including upper and lower bounds, least upper bounds, greatest lower bounds, bounded and unbounded functions, as well as simple proof and construction tasks.

The instructional structure follows an inquiry-oriented approach. Rather than introducing formal definitions immediately, students first explore concrete examples and identify patterns through guided questions. Formal definitions are introduced only after students have engaged with examples and simpler conceptual tasks. The system therefore aims to support conceptual understanding through progressively structured questioning. This approach is consistent with the findings of \citet{Bloch2003}, suggesting that concept formation based on concrete examples may support a deeper understanding of boundedness.

The application includes several question formats, including multiple-choice questions, ordering tasks, numerical input tasks, and short word-entry responses. In numerical and conceptual tasks, the application reacts differently to different types of incorrect answers and provides differentiated guiding questions.

\begin{figure}[H]
\centering
\includegraphics[width=.8\textwidth]{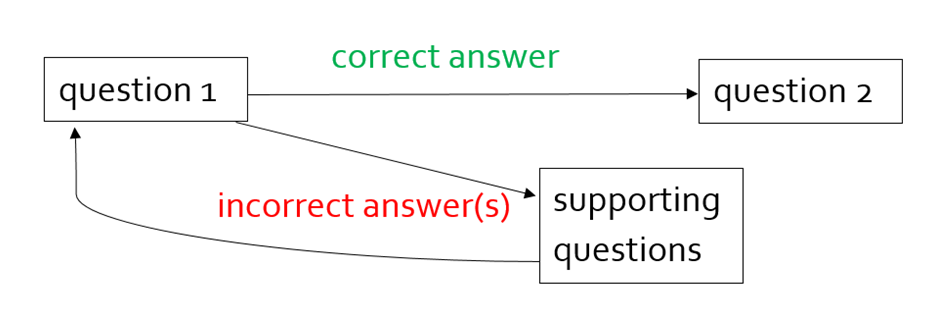}
\caption{Example interface of the autonomous inquiry-based learning application.}
\label{fig:interface}
\end{figure}

A central feature of the system is that it does not directly provide solutions when students answer incorrectly. Instead, it attempts to guide students toward the correct reasoning process through simpler intermediate questions. If a student demonstrates missing prerequisite knowledge, the application identifies the difficulty and recommends seeking additional support before continuing.

The final sections of the application include proof-based and construction tasks intended to support deeper conceptual understanding. These tasks were designed partly in response to difficulties identified in previous research, particularly students' confusion between boundedness and monotonicity and their difficulties interpreting unboundedness \citep{Bloch2003}.

Throughout the learning process, the application records students' responses, completion paths, and time spent on each task. These interaction logs were used to analyse students' progression through the application, while learning outcomes were evaluated separately using the paper-based assessment and compared with those of the control group.
\section{The Usefulness of the Application in Light of Advances in Artificial Intelligence}

When we began developing the application, artificial intelligence was not nearly as advanced as it is today. Since then, however, the question has naturally arisen whether manually designed educational systems of this type are still necessary, given the rapid progress of large language models and their increasing mathematical capabilities. We therefore considered it important to examine how current AI systems perform within the topic area addressed by the application.

To explore this question, we conducted several exploratory interactions with ChatGPT, using both logged-in and logged-out versions, in Hungarian and English. Overall, the system was often able to provide plausible explanations and solve routine problems; however, it also produced mathematically incorrect statements and internally inconsistent reasoning, sometimes maintaining false conclusions with high confidence.

For example, when discussing upper bounds of the set $\{3,5,7,9\}$, the model repeatedly identified 10 as the ``tightest'' upper bound and continued to defend this answer even after being prompted to compare it with smaller candidates such as 9.9. In another interaction, the model repeatedly claimed that no function defined on a closed interval could be unbounded, incorrectly applying the Weierstrass theorem even in situations where continuity was not assumed. Only after several additional prompts was it able to construct a valid counterexample.

These observations suggest that current LLM-based systems should be used cautiously as autonomous learning tools in mathematics education. While they are often capable of producing plausible explanations, they do not consistently support the guided conceptual discovery that inquiry-based learning requires.

The empirical results presented later in this paper further demonstrate that an inquiry-oriented questioning approach can effectively support students' conceptual understanding of boundedness. Together, these findings suggest that future AI-based educational systems may benefit from incorporating adaptive questioning structures that guide students through reasoning processes rather than primarily generating completed solutions.

\section{Methods}

\subsection{Participants}

A total of $N=92$ Hungarian secondary school students preparing for the advanced-level mathematics final examination participated in the study. Participants were recruited from six Budapest secondary schools whose advanced-level examination results were close to the national average. Students represented Grades 9--12 and were between 14 and 18 years old.

A matched control group consisting of 56 Hungarian secondary school students also participated in the assessment. These students were of the same age range, studied advanced-level mathematics, and attended another secondary school with a comparable educational profile. The control group received no instruction on boundedness before completing the assessment.

Participation took place in intact classes. An additional 42 students preparing for the intermediate-level mathematics examination also completed the application, although their results are not analysed further in the present paper.


\subsection{Procedure}

The study examined the effectiveness of the application in teaching the concept of boundedness of functions in situations where students had not previously encountered the topic during formal instruction. No pre-test was administered because teachers confirmed that boundedness had not yet been introduced in the participating classes, and students themselves also reported no prior instruction on the topic.

After a short introduction to the application and the generation of anonymous identification codes, students completed the learning material independently in a computer laboratory. They were given 90 minutes to work with the application.

The system automatically recorded

\begin{itemize}
\item students' responses,
\item completion paths,
\item time spent on each task,
\item progression through the adaptive question structure.
\end{itemize}

Some uncontrolled classroom factors, including varying levels of engagement and occasional peer interaction, may have influenced students' performance. Consequently, the findings should be interpreted conservatively.

Approximately one week after the intervention, the treatment group completed a paper-based follow-up assessment. During the intervening week, no additional instruction, homework, or revision related to boundedness was provided.

The control group completed exactly the same assessment under identical testing conditions using the same scoring rubric, but without having studied boundedness beforehand.


\subsection{Assessment}

To evaluate students' understanding of boundedness, six assessment tasks were developed according to the six cognitive levels of the revised Bloom taxonomy \citep{Anderson2001}. Since no standardized assessment instrument measuring conceptual understanding of boundedness was available, all tasks were specifically designed for this study.

The assessment measured students' ability to

\begin{itemize}
\item recall formal definitions,
\item classify functions,
\item apply concepts,
\item construct examples,
\item construct definitions,
\item produce mathematical arguments.
\end{itemize}

Table~\ref{tab:assessment} summarises the assessment tasks.

\begin{table}[H]
\centering
\caption{Assessment tasks aligned with the revised Bloom taxonomy.}
\label{tab:assessment}

\begin{tabular}{clp{8cm}}
\toprule
Task & Bloom level & Description\\
\midrule
1 & Remember &
Identify the definition of boundedness from below.\\

2 & Understand &
Classify functions according to boundedness.\\

3 & Apply &
Determine the greatest lower bound of a function.\\

4 & Create &
Construct a function satisfying prescribed conditions.\\

5 & Analyze &
Construct the definition of a function bounded from below.\\

6 & Evaluate / Create &
Prove that a function is unbounded from below.\\
\bottomrule
\end{tabular}

\end{table}

Tasks were assigned to Bloom levels according to the dominant cognitive process required for successful completion.

A Bloom level was regarded as achieved if the student successfully completed every task associated with the preceding levels. For Task~2, responses were accepted if at least three out of four classifications were correct.

Assessment was scored manually using an objective scoring rubric. No partial credit was awarded. However, computational mistakes that did not indicate conceptual misunderstanding were not counted against students when determining the highest Bloom level achieved.

In addition to the paper-based assessment, interaction data recorded automatically by the application—including completion time and progression through the adaptive learning sequence—were analysed statistically.

\section{Results and Discussion}

The effectiveness of the application was evaluated by analysing students' performance across the six levels of the revised Bloom taxonomy and by comparing the results with those of the matched control group.

\subsection{Learning outcomes across Bloom levels}

The first analysis examined the proportion of students reaching each cognitive level within the treatment group.

A Bloom level was considered achieved if the student correctly solved all tasks from the preceding levels. The only exception concerned Task~2, where students were considered successful if they correctly classified at least three of the four functions.

\begin{figure}[H]
\centering
\includegraphics[width=.85\textwidth]{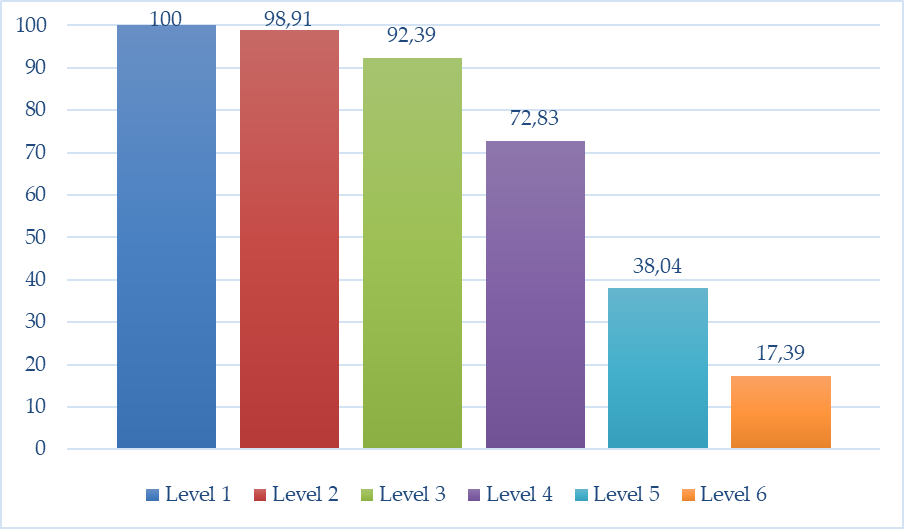}
\caption{Percentage of students achieving each Bloom level within the treatment group.}
\label{fig:bloomlevels}
\end{figure}

Figure~\ref{fig:bloomlevels} presents the percentage of students reaching each Bloom level.

To compare the treatment and control groups, two-sided Fisher's exact tests were performed for each Bloom level. Fisher's exact test was selected instead of the two-proportion z-test because several contingency tables contained small expected cell frequencies, including zero counts at the higher Bloom levels in the control group. Under such conditions, Fisher's exact test provides more reliable exact p-values than asymptotic methods.

\begin{table}[ht]
\centering
\caption{Comparison of treatment and control groups across Bloom levels using Fisher's exact test.}
\label{tab:fisher}
\begin{tabular}{cccc}
\toprule
\textbf{Bloom level} & \textbf{Treatment (\%)} & \textbf{Control (\%)} & \textbf{$p$ (Fisher's exact test)} \\
\midrule
1 & 100.00 & 69.64 & $1.16 \times 10^{-8}$ \\
2 & 98.91 & 53.57 & $2.13 \times 10^{-12}$ \\
3 & 92.39 & 7.14 & $2.90 \times 10^{-27}$ \\
4 & 72.83 & 1.79 & $9.42 \times 10^{-20}$ \\
5 & 38.04 & 0.00 & $3.18 \times 10^{-9}$ \\
6 & 17.39 & 0.00 & $5.29 \times 10^{-4}$ \\
\bottomrule
\end{tabular}
\end{table}
Students who completed the inquiry-based learning application substantially outperformed the matched control group at every Bloom level. Fisher's exact tests confirmed statistically significant differences between the treatment and control groups across all six cognitive levels.

The largest differences occurred at Bloom levels~3 and~4. More than 90\% of the treatment group successfully applied the concept of boundedness, whereas only 7.14\% of the control group reached this level. Similarly, nearly three quarters of the treatment group reached Bloom level~4, compared with fewer than 2\% of the control group.

The inclusion of a matched control group considerably strengthens the interpretation of these findings. Although participants were not randomly assigned and no pre-test was administered, the consistently large differences observed across every Bloom level suggest that the application contributed substantially to students' conceptual understanding.


\subsection{Completion of the application}

The second analysis investigated whether completing the application within the available time was associated with improved learning outcomes.

During classroom observations it became apparent that not all students remained equally engaged throughout the intervention. Some students occasionally interrupted their work or spent time on activities unrelated to the application. Consequently, completion status was examined separately.

A total of 70 students completed the entire inquiry sequence within the available time. Every participating school contributed at least one such student.

\begin{figure}[H]
\centering
\includegraphics[width=.85\textwidth]{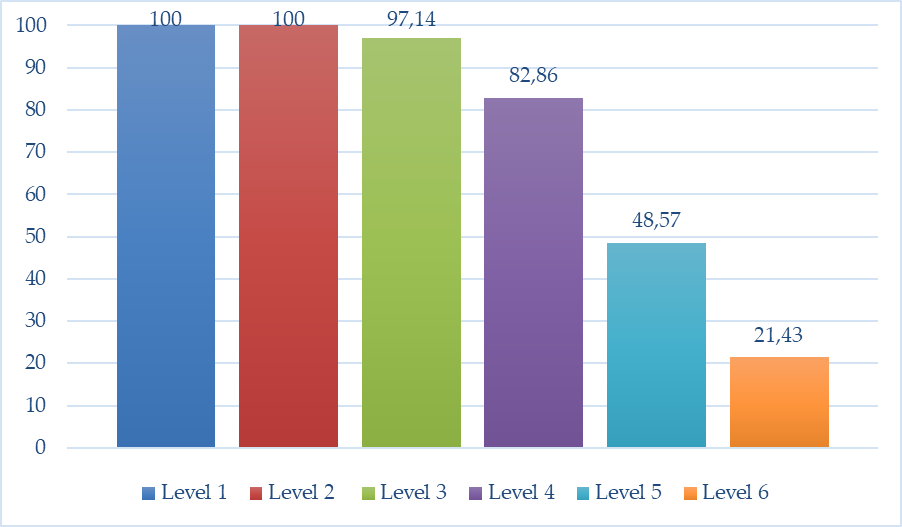}
\caption{Bloom levels achieved by students who completed the application.}
\label{fig:completed}
\end{figure}

Figure~\ref{fig:completed} shows the proportion of completers reaching each Bloom level.

\begin{figure}[H]
\centering
\includegraphics[width=.85\textwidth]{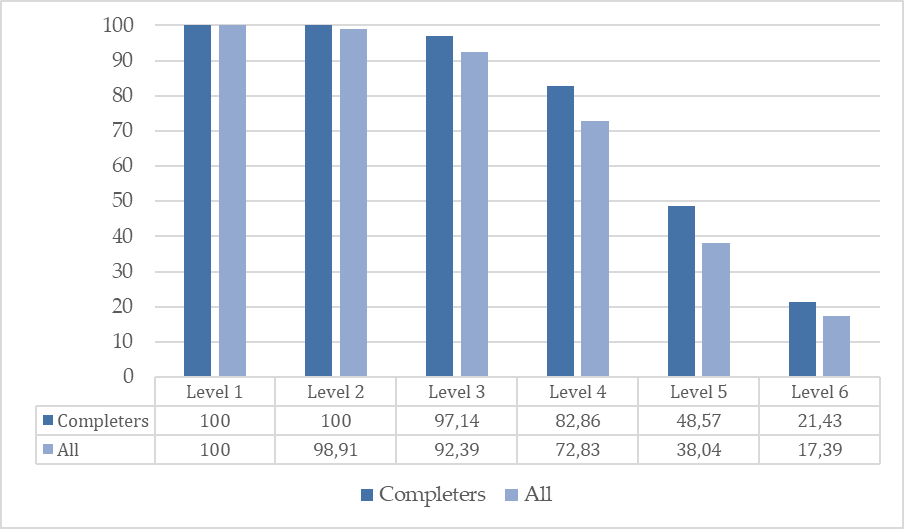}
\caption{Comparison between completers and the full sample.}
\label{fig:comparisoncompleted}
\end{figure}

For comparison, Figure~\ref{fig:comparisoncompleted} overlays the results of completers with those of the complete sample.

Students who completed the application generally achieved higher cognitive levels than the full sample.

Several explanations are possible. Students with stronger mathematical backgrounds may have progressed more rapidly through the adaptive sequence. At the same time, completing all inquiry tasks may itself have contributed to a deeper conceptual understanding of boundedness by exposing students to the full instructional sequence.

Regardless of the underlying mechanism, the results indicate that students who successfully completed the application almost always mastered the fundamental concepts and frequently achieved higher-order cognitive outcomes as well.


\subsection{Completion time}

The average time spent using the application was 1:12:55, while the median completion time was 1:19:58.

The shortest recorded completion time was 31:52, whereas the longest was 1:41:30.

To explore whether faster completion was associated with stronger performance, additional analyses were conducted for students finishing within

\begin{itemize}
\item 70 minutes,
\item 60 minutes,
\item 50 minutes.
\end{itemize}

These analyses should be interpreted cautiously because the corresponding groups were considerably smaller ($n=36$, $27$, and $12$, respectively) and represented fewer participating schools.

Nevertheless, some consistent tendencies emerged.

\begin{figure}[H]
\centering
\includegraphics[width=.90\textwidth]{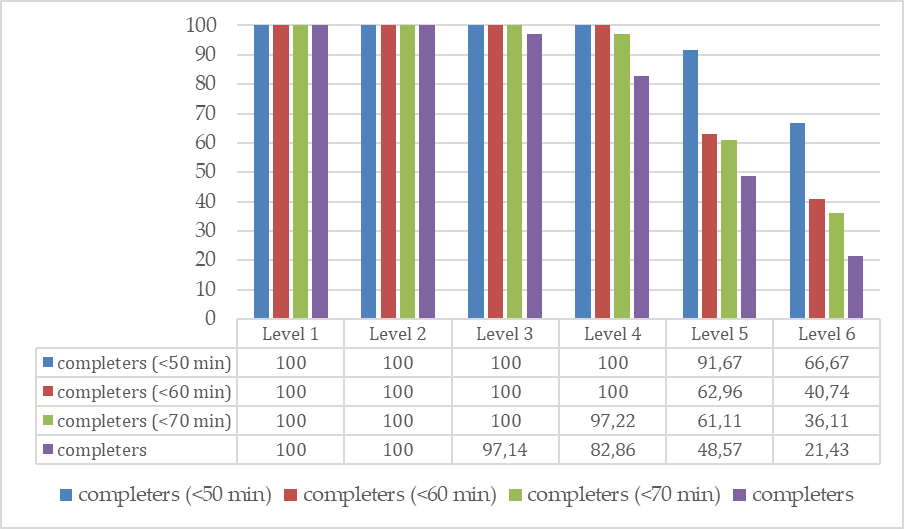}
\caption{Comparison of learning outcomes for students completing the application within different time limits.}
\label{fig:timecomparison}
\end{figure}

Figure~\ref{fig:timecomparison} compares the results of these groups with those of all students completing the application.

For the lower Bloom levels, shorter completion times were generally associated with a higher proportion of successful students.

Within the present sample, every student completing the application within 60 minutes reached at least Bloom level~4.

Although the sample sizes are insufficient for strong statistical conclusions, these findings raise an interesting possibility that the application could eventually support not only learning but also aspects of formative assessment. Future research may therefore investigate whether completion behaviour and interaction patterns could provide useful indicators of students' conceptual understanding.

\section{Conclusion}

This study investigated the effectiveness of an autonomous inquiry-oriented learning application designed to teach the concept of boundedness of functions to secondary school students. The application was developed to support situations in which teacher guidance is unavailable or limited by replacing direct explanation with adaptive questioning and structured conceptual guidance.

The first research question examined whether the application improved students' understanding of boundedness compared with students who had not received instruction on the topic and what levels of understanding students were able to achieve. The results provide strong evidence that students using the application substantially outperformed a matched control group across every level of the revised Bloom taxonomy. Nearly all participants mastered the fundamental concepts, most successfully applied them in unfamiliar situations, and a considerable proportion also demonstrated higher-order reasoning through analysis, proof, and construction tasks. Although the study did not employ random assignment or a pre-test, the magnitude and consistency of the observed differences suggest that the application made a substantial contribution to students' conceptual understanding.

The second research question addressed the relationship between students' performance and the pace at which they progressed through the application. Students who completed the entire inquiry sequence within the available time generally achieved higher cognitive levels on the follow-up assessment than the overall sample. While this relationship cannot be interpreted causally, it suggests that successfully progressing through the complete sequence of adaptive questions may be associated with deeper conceptual understanding. At the same time, faster completion may also reflect stronger prior mathematical preparation or greater engagement during the learning process. Future studies incorporating pre-tests and additional learner characteristics could help disentangle these explanations.

Beyond evaluating a specific instructional application, the study contributes to the broader discussion concerning the role of artificial intelligence in mathematics education. Contemporary large language models have demonstrated impressive mathematical capabilities, yet they are primarily designed to generate answers rather than to guide students through carefully structured reasoning processes. Our exploratory observations suggest that, despite their strengths, current systems do not consistently provide the kind of adaptive conceptual guidance required for inquiry-based learning. In contrast, the application investigated here deliberately avoids presenting immediate solutions and instead supports students through progressively structured questions and differentiated feedback.

These findings suggest that future AI-supported educational systems may benefit from integrating inquiry-oriented instructional principles into modern artificial intelligence. Rather than viewing adaptive tutoring systems and large language models as competing technologies, future educational environments may combine the strengths of both approaches. Large language models could provide flexible natural-language interaction, whereas inquiry-oriented pedagogical frameworks could regulate the learning process by determining which questions should be asked, when additional guidance should be provided, and how conceptual understanding should gradually be constructed. Such hybrid systems may offer a promising direction for future research in artificial intelligence for mathematics education.

Several limitations should be acknowledged. Participants were recruited from intact classes rather than randomly assigned to experimental conditions, and no pre-test was administered. Although the matched control group substantially strengthens the interpretation of the findings, causal conclusions should nevertheless be drawn with appropriate caution. In addition, the study focused exclusively on a single mathematical topic and measured learning approximately one week after the intervention. Further research should therefore examine long-term retention, transfer to other mathematical concepts, and the effectiveness of similar inquiry-oriented applications across different areas of mathematics and in more diverse educational settings.

Overall, the findings indicate that autonomous inquiry-based learning environments can effectively support meaningful mathematical learning without continuous teacher intervention. More broadly, they suggest that the future development of artificial intelligence in mathematics education may benefit not only from increasingly powerful answer-generation capabilities but also from instructional designs that encourage students to think, explore, and construct mathematical knowledge through carefully guided inquiry.

\section*{Declaration of Generative AI Use}

ChatGPT (OpenAI) was used to improve the English language and readability of the manuscript. 
The author reviewed and edited all AI-assisted suggestions and takes full responsibility for the final content of the manuscript.
\bibliographystyle{plainnat}
\bibliography{references}

\end{document}